\documentclass[12pt]{article}
\usepackage{amsmath}
\usepackage{latexsym}
\usepackage{amssymb}
\newtheorem{thm}{Theorem}[section]
\newtheorem{la}[thm]{Lemma}
\newtheorem{Defn}[thm]{Definition}
\newtheorem{Remark}[thm]{Remark}
\newtheorem{prop}[thm]{Proposition}

\newtheorem{Example}[thm]{Example}
\newtheorem{Examples}[thm]{Examples}
\newtheorem{Number}[thm]{\!\!}
\newenvironment{defn}{\begin{Defn}\rm}{\end{Defn}}

\newenvironment{rem}{\begin{Remark}\rm}{\end{Remark}}
\newenvironment{numba}{\begin{Number}\rm}{\end{Number}}
\newenvironment{proof}{{\noindent\bf Proof.}}%
                  {\nopagebreak\hspace*{\fill}$\Box$\medskip\medskip\par}   
\newcommand{\Punkt}{\nopagebreak\hspace*{\fill}$\Box$}
\newcommand{\wb}{\overline}
\newcommand{\ve}{\varepsilon}

\newcommand{\isom}{\cong}
\newcommand{\mto}{\mapsto}
\newcommand{\N}{{\mathbb N}}

\newcommand{\R}{{\mathbb R}}

\newcommand{\C}{{\mathbb C}}

\newcommand{\cB}{{\mathcal B}}

\newcommand{\cL}{{\mathcal L}}

\newcommand{\cg}{{\mathfrak g}}
\newcommand{\ch}{{\mathfrak h}}

\newcommand{\one}{{\bf 1}}
\newcommand{\sub}{\subseteq}
\DeclareMathOperator{\GL}{GL}
\DeclareMathOperator{\im}{im}

\DeclareMathOperator{\id}{id}

\newcommand{\sbull}{{\scriptscriptstyle \bullet}}

\DeclareMathOperator{\Diff}{Diff}

\DeclareMathOperator{\Hom}{Hom}
\DeclareMathOperator{\Gf}{Gf}
\DeclareMathOperator{\Gh}{Gh}
\newcommand{\dl}{{\displaystyle\lim_{\longrightarrow}}}

\begin{document}
\begin{center}
{\Large\bf Solutions to
open problems in Neeb's recent\\[1.6mm]
survey on infinite-dimensional Lie groups}\\[5mm]
\renewcommand{\thefootnote}{\fnsymbol{footnote}}
{\bf Helge Gl\"{o}ckner\footnote{\,Supported by the German Research Foundation
(DFG), projects GL 357/5-1
and GL~357/7-1}}\vspace{4mm}
\end{center}
\renewcommand{\thefootnote}{\arabic{footnote}}
\setcounter{footnote}{0}
\begin{abstract}
\hspace*{-5.8mm}We solve three open problems concerning
infinite-dimensional Lie\linebreak
groups posed in
a recent survey article
by K.-H. Neeb~\cite{Jap}. Moreover,
we prove
a result by the author
announced in~\cite{Jap},
which answers a
question posed
in an earlier, unpublished  version of the survey.\vspace{.8mm}
\end{abstract}
{\footnotesize {\em Classification}:
22E65 (Primary) 58B10 (Secondary)\\[2.5mm]
{\em Key words}:
infinite-dimensional Lie group, open problem,
initial submanifold, initial Lie subgroup,
direct limit group, topological group with Lie algebra,
Borel's Theorem, Taylor series,
diffeomorphism}
\begin{center}
{\bf\Large Introduction}
\end{center}
The recent article~\cite{Jap}
by K.-H. Neeb surveyed the state of the art
of the\linebreak
theory of infinite-dimensional Lie groups,
and described typical open\linebreak
problems.
In this article, we provide full solutions to three of these:
\begin{itemize}
\item
We describe a subgroup of an infinite-dimensional
Lie group which is not an initial Lie subgroup (Section~\ref{secnon}).
This answers
\cite[Problem II.6]{Jap},
the ``Initial Subgroup Problem,''
in the negative.
\item
We show that every subgroup of a direct limit
$G=\dl\, G_n$\vspace{-.7mm}
of finite-dimensional
Lie groups is an initial Lie subgroup
(Section~\ref{yesinit}).
This answers \cite[Problem VII.3]{Jap}
in the affirmative.
\item
We show that every direct limit
$G=\dl\, G_n$\vspace{-.7mm}
of finite-dimensional Lie groups
is a topological group with Lie algebra
in the sense of~\cite{HaM}\linebreak
(Section~\ref{secgpLiea}).
This answers \cite[Problem VII.2]{Jap}
in the affirmative.
\end{itemize}
We also prove a result by the author
announced in \cite[Proposition~VI.2.11]{Jap}:
Let $M$ be $\sigma$-compact finite-dimensional
smooth manifold
of dimension~$n>0$,
and $\Diff_c(M)$
be the Lie group of smooth diffeomorphisms $\gamma\colon M\to M$
which are compactly supported
in the sense that the closure of
$\{x\in M\colon \gamma(x)\not=x\}$
is compact (see \cite{Mic} and \cite{DFM}).
Given $p\in M$, let
$\Diff_c(M)_p$ be its stabilizer,
and $\Diff_c(M)_{p,0}$
be the subgroup of all
$\gamma\in \Diff_c(M)_p$ such that
$\det T_p(\gamma)>0$,
where $T_p(\gamma)\colon \hspace*{-.8mm}
T_p M \hspace{-1mm}\to\hspace{-1mm} T_p M$
is the tangent map.
Let
$\Gf_n(\R)\sub \R[\![x_1,\ldots, x_n]\!]$
be the
group of all formal diffeomorphisms
of $\R^n$,
and $\Gf_n(\R)_0$
be the subgroup of all
formal diffeomorphisms whose
linear part has positive
determinant (see \cite[pp.\ 426--427]{Jap} and \cite{RaK}).
In Section~\ref{secborel},
we establish
a version of \'{E}. Borel's theorem
on the existence of smooth maps
with a given Taylor series:\\[2.5mm]
{\bf Borel's Theorem for Diffeomorphisms.}
\emph{Let a finite-dimensional smooth manifold
$M$ of dimension $n$ be given,
$p\in M$
and $\phi$ be a local chart around~$p$
such that $\phi(p)=0$.
For each $\gamma \in \Diff_c(M)$,
let $T^\infty_p(\gamma)$ denote the Taylor series of $\phi\circ \gamma\circ \phi^{-1}$
around~$0$. Then}\vspace{-1.5mm}
\[
T^\infty_p\colon \Diff_c(M)_{p,0}\to \Gf_n(\R)_0\vspace{-1.5mm}
\]
\emph{is a surjective homomorphism of Lie groups.}\\[2.5mm]
We mention that further problems from \cite{Jap}
have an affirmative answer:\\[2.5mm]
In \cite{DaG}, it is shown that the group
$\Gh_n(\C)$ of germs of analytic
diffeomorphisms around $0\in \C^n$
(as first considered in~\cite{Pis})
is a regular Lie group in Milnor's sense.\footnote{Likewise
for groups of germs
of diffeomorphisms around a compact set, as~in~\cite{JFA}.}
This answers the first half of
\cite[Problem~VI.5]{Jap} in the affirmative.\\[2.5mm]
Using a version of the Frobenius Theorem
for co-Banach distributions on infinite-dimensional
manifolds,\footnote{The adaptation to Keller's $C^\infty_c$-theory
(our setting of infinite-dimensional calculus)
of Hiltunen's Frobenius Theorem from~\cite{Hil}
(to be found in~\cite{GaL}).}
one obtains the following result
concerning integral subgroups in the sense
of \cite[Definition~IV.4.7]{Jap}
(see \cite{GaL}):\\[2.5mm]
\emph{Let $G$ be a Lie group modelled on a locally
convex space and
$\ch\sub L(G)$ be a Lie subalgebra
such that $L(G)=\ch\oplus X$
for some Banach space $X\sub L(G)$.
Then there exists an integral subgroup
$H\sub G$ such that $L(H)=\ch$.}\\[2.5mm]
Note that a complement $X$ always exists if $\ch$ has finite
codimension. Hence \cite[Problem~VI.4]{Jap}
has an affirmative answer (even if $G$
is not regular).\\[2.5mm]
Finally, we mention that the idea expressed in
\cite[Problem~VII.5]{Jap}
has been taken up in
\cite[Lemma~14.5 and Section~15]{JFA}.

\noindent
{\bf Prerequisites.}
The reader is referred to \cite{Jap},
\cite{GaN} or \cite{RES}
for the setting of
infinite-dimensional Lie groups, manifolds
and $C^k$-mappings used in this\linebreak
article
(cf.\ \cite{Mil} for the case
of sequentially complete modelling spaces).
The words ``Lie group'' and ``manifold''
will always refer to Lie groups and\linebreak
manifolds
modelled on locally convex spaces.
\section{Example
of a subgroup of a Lie group
which is not an initial Lie subgroup}\label{secnon}
We shall use the following terminology.
\begin{defn}\label{definit}
Let $\N:=\{1,2,\ldots\}$, $k\in \N \cup\{\infty\}$
and $M$ be a $C^k$-manifold.
A subset $N\sub M$ is called a \emph{$C^k$-initial
submanifold} of~$M$
if $N$ can be made a $C^k$-manifold
with the following properties:
\begin{itemize}
\item[(a)]
The inclusion map
$\lambda \colon N\to M$ is a $C^k$-map; and
\item[(b)]
If $X$ is a $C^k$-manifold
and $f\colon X\to N$ a map
such that $\lambda \circ f$ is $C^k$,
then~$f$ is~$C^k$.\vspace{.9mm}
\end{itemize}
Such a manifold structure is called
\emph{$C^k$-initial} in~$M$.
If $G$ is a Lie group, let $L(G):=T_1(G)$
be its Lie algebra.
Following \cite[Definition~II.6.1]{Jap},
we say that a subgroup
$H\sub G$ is an \emph{initial Lie subgroup}
if it can be given a smooth manifold structure
which makes $H$ a Lie group,
is $C^k$-initial in~$G$
for each $k\in \N\cup\{\infty\}$,
and gives rise to an injective
tangent map $L(\lambda):=T_1(\lambda)\colon
L(H)\to L(G)$.
\end{defn}
\begin{rem}\label{betterinit}
It means additional information
if $N$ and $X$ may even be $C^k$-manifolds
with $C^k$-boundary in Definition~\ref{definit},
or more generally
$C^k$-manifolds ``with rough boundary''
as in~\cite{GaN}.\footnote{The latter manifolds
locally look like
a locally convex subset with dense
interior in the modelling space.
They include manifolds with $C^k$-boundary
and manifolds with corners
as special cases.}
We then call $N$ a \emph{strongly initial
$C^k$-submanifold} of~$M$.
Notably, the strengthened concept
enables maps $f\colon X\to N$
to be considered
if $X=[0,1]$, $X=[0,1]^2$
or $X$ is a simplex.
\end{rem}
Let $\R^\N$ be the Fr\'{e}chet space of
real sequences (with the product topology).
We show:\vspace{1mm}
\begin{thm}\label{thmnoin}
For each $k\in \N
\cup\{\infty\}$,
the additive subgroup $N:=\ell^\infty$
of bounded sequences is not a $C^k$-initial
submanifold of $M:=\R^\N$.
In particular, it is not an initial Lie subgroup.
\end{thm}
\begin{proof}
We shall use a well-known fact:
A map from a $C^k$-manifold
to $\R^\N$ (or any direct product
of locally convex spaces) is $C^k$ if and only if
all of its components are~$C^k$
(see \cite[Chapter~1]{GaN}
or \cite[Lemma~10.2]{Ber}).\\[2.5mm]
It is useful to note first that,
for each $m\in \N$,
the topological product
$K_m:=[{-m},m]^\N$
is a metrizable
compact topological space
(and therefore\linebreak
separable),
whence every open
subset $U\sub K_m$ is $\sigma$-compact.\\[2.5mm]
If the theorem was false,
we could equip $N$ with a $C^k$-initial
manifold structure.
Assuming this,
write $\lambda\colon N\to M$ for
the inclusion map (which is $C^k$).
Given $m\in \N$,
choose an injective $C^k$-map $h_m\colon \R\to \;]{-2m},2m[$
such that $h_m([{-1},1])=[{-m},m]$.
Then the injective map
\[
g_m\colon \R^\N\to N\,,\quad (x_n)_{n\in \N}\mto
(h_m(x_n))_{n\in \N}
\]
is $C^k$ by initiality,
exploiting that $\lambda\circ g_m$
is $C^k$ since its components are~$C^k$.
In particular, $g_m$ is continuous,
entailing that
\[
g_m|_{[{-1},1]^\N}\colon [{-1},1]^\N\to N
\]
is a topological embedding.
Since $g_m([{-1},1]^\N)=K_m$
as a set and $\lambda\circ g_m$
takes $[{-1},1]^\N$ homeomorphically onto
$K_m$, it follows that $N$ induces
the given compact topology on~$K_m$.\\[2.5mm]
By the preceding, $N=\bigcup_{m\in \N}K_m$
is $\sigma$-compact.
Let $E$ be the modelling
space of~$N$ and $\phi\colon N\supseteq U\to V\sub E$
be a chart of~$N$, defined on a non-empty open set $U\sub N$.
Then $U=\bigcup_{m\in \N}(U\cap K_m)$
where $U\cap K_m$ is $\sigma$-compact
for each $m\in \N$ (as recalled above),
and hence $U$ is $\sigma$-compact.
Consequently, also the locally convex
space~$E$ is $\sigma$-compact,
because it contains a $\sigma$-compact set
with non-empty interior.\\[2.5mm]
If $x\in N$, then $T_x(N)\isom E$
is $\sigma$-compact.
We claim that the tangent map\linebreak
$T_x(\lambda)\colon T_x(N)\to T_x(M)$
is surjective.\\[2mm]
If this is true, then $\R^\N\isom T_x(M)=T_x(\lambda)(T_x(N))$
is $\sigma$-compact. But this is absurd
(by a Baire argument),
and the theorem is proved
by contradiction.
To prove the claim, identify $T_x(M)$ with $\R^\N$.
If $y=(y_n)_{n\in \N}\in \R^\N$,
then
\[
g \colon \R\to N\,,\quad g(t)\, :=\, (x_n+\sin(t y_n ))_{n\in \N}
\]
is a function such that $\lambda\circ g\colon \R\to \R^\N$
is $C^k$, because
all components of this composition are $C^k$.
Furthermore, $T_x(\lambda)(g'(0))=(\lambda\circ g)'(0)
=y$.
\end{proof}
Our proof
shows that Theorem~\ref{thmnoin}
remains valid if, instead of our general
setting, only
manifolds modelled on complete
(or sequentially complete,
or Mackey complete) locally
convex spaces are considered.
Also, it remains valid if manifolds
are replaced by
manifolds with rough boundary.
\section{Subgroups of direct limit groups
are initial}\label{yesinit}
If $G_1\sub G_2\sub \cdots$
is a sequence of finite-dimensional real
Lie groups such that the inclusion maps
$G_n\to G_{n+1}$ are smooth homomorphisms,
then $G:=\bigcup_{n\in \N}G_n$
admits a Lie group structure
making each inclusion map $G_n\to G$ a smooth homomorphism,
and such that $G=\dl\, G_n$\vspace{-.8mm}
in the category of smooth Lie groups
modelled on locally convex spaces \cite[Theorem~4.3]{FUN}
(cf.\ \cite{NRW},
\cite[Theorem~47.9]{KaM} and \cite{KYO} for special cases).
Then $G=\dl\,G_n$\vspace{-.8mm}
also as a topological space,
and as a $C^k$-manifold, for each $k\in \N_0\cup\{\infty\}$.
Furthermore,\vspace{-1mm}
\[
L(G)\; =\;  {\bigcup}_{n\in \N}\, L(G_n) \; =\; \dl\,L(G_n)\vspace{-1mm}
\]
as a locally convex space,
topological space and $C^k$-manifold
(see\linebreak
\cite[Theorems~4.3\,(a) and 3.1]{FUN} for all of this).
For brevity, we shall refer to
$G=\dl\,G_n$\vspace{-1mm}
as a \emph{direct limit Lie group}.
In this section, we show:
\begin{thm}\label{mainXX}
Let $G_1\sub G_2\sub \cdots$ be finite-dimensional Lie groups,
such that the inclusion map $G_n\to G_{n+1}$
is a smooth homomorphism
for each $n\in \N$.
Let $G=\dl\,G_n$\vspace{-.8mm}
be the direct limit Lie group,
and $H\sub G$ be a subgroup.
Equip $H_n:=G_n\cap H$ with
the finite-dimensional real
Lie group structure induced by~$G_n$
$($as in {\rm \cite[Chapter~III, \S4.5, Proposition~9]{Bou})}
and $H$ with the Lie group
structure making it the direct limit
Lie group $H=\dl\, H_n$\vspace{-.8mm}
$($as in
{\rm\cite[Proposition~7.2]{FUN}).}
Then the manifold structure on~$H$
is strongly $C^k$-initial in~$G$, for all
$k\in \N\cup\{\infty\}$.
Thus $H$ is an initial Lie subgroup
of~$G$.
\end{thm}
The final conclusion is possible because
the inclusion map $\lambda\colon H\to G$
gives rise to an embedding
$L(\lambda) \colon L(H)\to L(G)$
by \cite[Proposition~7.2]{FUN}.\\[2.5mm]
It is useful to recall some simple facts.
\begin{numba}\label{fct1}
The identity component
of a direct limit Lie group $G=\dl\, G_n$\vspace{-.8mm} is
\[
G_0\; =\; \bigcup_{n\in \N}(G_n)_0\; =\; \dl\,(G_n)_0\,.
\]
In particular, $G_0$ is a direct
limit group as well.
\end{numba}
\begin{numba}\label{fct2}
Consider an
ascending sequence
$X_1\sub X_2\sub\cdots$
of topological spaces.
Then
$X:=\bigcup_{n\in \N}X_n$
is its direct limit in the category of
topological spaces, if we declare $U\sub X$
open if and only if $U\cap X_n$ is open in~$X_n$
for each $n\in \N$. The topology so obtained
is called the \emph{direct limit topology}.
If $U_n\sub X_n$ are open subsets
such that $U_1\sub U_2\sub\cdots$,
then $U:=\bigcup_{n\in \N} U_n$
is open in~$X$,
because $U\cap X_n=\bigcup_{m\geq n}U_m\cap X_n$
is open in $X_n$ for each~$n$.
\end{numba}
Also the following lemma
will be useful for the proof of Theorem~\ref{mainXX}.
\begin{la}\label{tooliniti}
Let $X$ be a
finite-dimensional real vector space,
$K\sub X$ be a non-empty compact set,
$k\in \N\cup\{\infty\}$, $M$ be a finite-dimensional
$C^k$-manifold
and $f\colon U\to M$ be a $C^k$-map
on an open neighbourhood $U$ of~$K$ in~$X$
such that $f|_K$ is injective and the tangent map
$T_x(f)\colon T_x(X)\to T_{f(x)}(M)$
is invertible for each $x\in K$.
Then there is an open neighbourhood $V\sub U$ of~$K$ such that $f(V)$
is open in~$M$ and $f|_V$ is a $C^k$-diffeomorphism
onto $f(V)$.
\end{la}
\begin{proof}
The proof of \cite[Lemma~15.6]{JFA}
(where $M=X$)
carries over.\footnote{Unfortunately,
there is a misprint:
``$E$'' reads~``$X$'' in the lemma.}
\end{proof}
{\bf Proof of Theorem~\ref{mainXX}.}
We may assume that $G_1=H_1=\{1\}$.
Let $X$ be a $C^k$-manifold (possibly with rough boundary)
and $f\colon X\to H$ be a map such that
$\lambda\circ f\colon X \to G$ is $C^k$,
where $\lambda\colon H\to G$ is the inclusion map.
Given $p\in X$, let $Y$ be the connected
component of~$X$ which contains~$p$.
Since $Y$ is open in~$X$,
it suffices to show that
$f(p)^{-1}f|_Y$ is $C^k$
on some open neighbourhood of~$p$.
After replacing~$X$ with $Y$ and $f$ with
$f(p)^{-1}f|_Y$, we may assume that $X$ is
connected and $f(p)=1$.
Then $f(X)$ is contained in the identity
component $G_0$
of~$G$.
Since $G_0=\dl\,(G_n)_0$\vspace{-.8mm}
(see {\bf\ref{fct1}}),
after replacing $G$ with $G_0$,
$G_n$ with $(G_n)_0$ and $H$ with $H\cap G_0$,
we may assume that $G$ and each $G_n$ is connected.
Then $G_n$ is the integral subgroup
of $G_{n+1}$ with Lie algebra $L(G_n)$.
Since the latter is $C^\infty$-initial
in $G_{n+1}$, it follows that
\begin{equation}\label{givesgoodH}
L(G_n)\; =\; \{x\in L(G_{n+1})\colon \exp_{G_{n+1}}(\R x)\sub G_n\}
\end{equation}
(see \cite[Theorem~IV.4.14]{Jap}).
We recall that
\begin{equation}\label{secondform}
L(H_n)\;=\; \{\gamma'(0)\colon \gamma \in C^1([0,1],G_n)
\;\mbox{with $\gamma(0)=1$ and $\gamma([0,1])\sub H_n$}\}
\end{equation}
and also
\begin{equation}\label{firstform}
L(H_n)\;=\;\{x\in L(G_n)\colon \exp_{G_n}(\R x)\sub H_n\}\,,
\end{equation}
(cf.\ \cite[Chapter~3, \S4.5]{Bou}
or \cite[Theorem~IV.4.14]{Jap}).
Combining (\ref{firstform})
with (\ref{givesgoodH}) and the fact
that $L(G)=\bigcup_{n\in \N}L(G_n)$,
we see that
\[
L(H_n)=\{x\in L(G)\colon \exp_G(\R x)\sub H_n\}
\]
and thus
\begin{equation}\label{veryuse}
L(H)\cap L(G_n)\; =\; L(H_n)\,.
\end{equation}
Recursively, we find vector subspaces
$E_n\sub L(G_n)$ for $n\in \N$
such that $E_{n-1}\sub E_n$ (if $n\geq 2$)
and
\[
L(G_n)=L(H_n)\oplus E_n\,.
\]
In fact, if $n=1$ we choose $E_1:=\{0\}$.
If $E_1,\ldots, E_n$ have
been chosen, then 
\[
L(H_{n+1})\cap E_n
= L(H_{n+1})\cap L(G_n)\cap E_n
= L(H_n)\cap E_n=\{0\}\,,
\]
using (\ref{veryuse}) for the penultimate
equality. Hence $E_n$ can be extended
to a vector complement $E_{n+1}$
to $L(H_{n+1})$ in $L(G_{n+1})$.\\[2.5mm]
Give $E:=\bigcup_{n\in \N}E_n=\dl\, E_n$\vspace{-.8mm}
the locally convex direct limit topology.
Then
\[
L(G)\; =\; \dl\, (L(H_n)\oplus E_n)
\; =\; \dl\, L(H_n)\oplus \dl\, E_n
\;=\; L(H)\oplus E\vspace{-.8mm}
\]
as a locally convex space (see, e.g., \cite[Theorem~3.4]{HST}).
For each $n\in \N$, let $R_n\sub L(H_{n+1})$
and $S_n\sub E_{n+1}$
be vector subspaces such that
\[
L(H_{n+1})\;=\; L(H_n)\oplus R_n\quad\mbox{and}\quad
E_{n+1}\; =\; E_n\oplus S_n\,.
\]
We now pick a basis $\cB$ of $L(G)$
such that $\cB\cap L(H_n)$
is a basis of $L(H_n)$,
$\cB\cap E_n$ is a basis of $E_n$,
$\cB\cap R_n$ is a basis of $R_n$
and $\cB\cap S_n$ is a basis of $S_n$,
for each $n\in \N$. It can be used to
introduce the supremum norm
\[
\|.\|\colon L(G)\to [0,\infty[\,,\quad \left\| \sum_{b\in \cB } t_b b\right\|
\; :=\; \sup_{b\in \cB} \, |t_b|
\quad\mbox{for $\; (t_b)_{b\in \cB}\in \R^{(\cB)}$.}
\]
Given $t>0$ and $n\in \N$, abbreviate
\[
U^n_t\;:=\; \{x\in L(H_n)\colon \|x\|<t\}\,,\qquad
V^n_t\; :=\; \{x\in E_n \colon \|x\|<t\}\,,
\]
\[
X^n_t\; :=\;  \{ x \in R_n \colon \|x\|<t\}\,,
\quad \mbox{and}\quad
Y^n_t\; :=\;  \{ x \in S_n \colon \|x\|<t\}\,.
\]
Then
$U^{n+1}_t=U^n_t\times X^n_t$ and
$V^{n+1}_t=V^n_t\times Y^n_t$.
Let $2=t_1>t_2>\cdots>1$.\\[2.5mm]
{\bf Claim.} \!\emph{There are $C^\infty$-maps $\eta_n\colon \! V^n_{t_n} \to G_n$
for $n\!\in \!\N$
and $C^\infty$-diffeomorphisms
$\gamma_n\colon U^n_{t_n} \to H_n$
onto open identity neighbourhoods in~$H_n$,
such that $\gamma_n(0)=\eta_n(0)=1$,}
\[
g_n\colon U^n_{t_n}\times V^n_{t_n}
\to G_n\,,\quad g_n(x,y):=\gamma_n(x)\eta_n(y)
\]
\emph{is a $C^\infty$-diffeomorphism onto an
open identity neighbourhood in $G_n$, and}
\begin{equation}\label{hconbin}
\gamma_k|_{U^\ell_{t_k}}\,=\,\gamma_\ell|_{U^\ell_{t_k}}\quad
\mbox{as well as}\quad
\eta_k|_{V^\ell_{t_k}}\,=\,\gamma_\ell|_{V^\ell_{t_k}}\quad
\mbox{for all $1\leq \ell\leq k\leq n$.}
\end{equation}
If this claim is true, then
$\sigma_n:=\gamma_n|_{U^n_1}$ is a $C^\infty$-diffeomorphism
onto the open identity neighbourhood $Q_n:=\gamma_n(U^n_1)$ in~$H_n$,
the mapping $h_n:=g_n|_{U^n_1\times V^n_1}$ is a $C^\infty$-diffeomorphism
onto the open identity neighbourhood~$g_n(U^n_1\times V^n_1)$\linebreak
$=:P_n$
in~$G_n$,
and $\tau_n:= \eta_n|_{V^n_1}\colon V^n_1\to G_n$
is a smooth map (and in fact an immersion
and topological embedding, since
$h_n(x,y)=\sigma(x)\tau(y)$ and~$h_n$ is a diffeomorphism.).
Then $U:=\bigcup_{n\in \N}U^n_1$ is an open $0$-neighbourhood
in $L(H)$
and
$Q:=\bigcup_{n\in \N}Q_n$ is a connected, open identity neighbourhood
in~$H$\linebreak
(see {\bf\ref{fct2}}). Moreover,
\[
\sigma:=\dl\, \sigma_n\colon L(H)\supseteq U\to Q\sub H\vspace{-.8mm}
\]
is a $C^\infty$-diffeomorphism (cf.\
\cite[Lemma~1.9 and Proposition~3.3]{FUN}).
Also,
$W:=\bigcup_{n\in \N}(U^n_1\times V^n_1)$ is an open $0$-neighbourhood
in $L(G)$,
$P:=\bigcup_{n\in\N}P_n$ is an open identity neighbourhood in~$G$,
and
\[
h:=\dl\, h_n\colon L(G)\supseteq W\to P\sub G\vspace{-.8mm}
\]
is a $C^\infty$-diffeomorphism.
Finally, $V:=\bigcup_{n\in \N} V_n$ is an open $0$-neighbourhood
in~$E$ and
\[
\tau:=\dl\, \tau_n\colon E\supseteq V\to G
\]
is a $C^\infty$-map (and actually an immersion
and topological embedding, since
$h(x,y)=\sigma(x)\tau(y)$ and~$h$ is a diffeomorphism).
There exists a connected open neighbourhood
$Z\sub X$ of~$p$ such that $f(Z)\sub P$.
We now show that
\begin{equation}\label{ZsubQ}
f(Z)\; \sub \; Q\,.
\end{equation}
To this end, let
$q\in Z$.
There exists a $C^1$-map
$\theta\colon [0,1]\to Z$
such that $\theta(0)=p$
and $\theta(1)=q$.
By \cite[Lemma~1.7\,(d)]{FUN}, there exists $n\in \N$ such that
\[
f(\theta([0,1]))\; \sub \; G_n\,.
\]
Then the left logarithmic
derivative $\delta(f\circ\theta)\colon [0,1]\to L(G)$
is a continuous map
taking its values in~$L(H_n)$
(by (\ref{secondform})
and \cite[Chapter~3, \S4.5, Lemma~4]{Bou}).\footnote{Recall
that if $\gamma\colon [0,1]\to G$
is $C^1$ then $\delta(\gamma)(t):=\gamma(t)^{-1}\cdot \gamma'(t)\in L(G)$
defines $\delta(\gamma)\colon [0,1]\to L(G)$
(cf.\ \cite[p.\,1043]{Mil} and \cite{Jap}).
The dot denotes multiplication in the tangent
group~$TG$.}
As a consequence of the standard Existence
and Uniqueness Theorem for solutions
to ordinary differential
equations,
there exists a unique $C^1$-curve $\kappa\colon
[0,1]\to (H_n)_0$
such that $\delta(\kappa)=\delta(f\circ \theta)$
and $\kappa(0)=1$.
Then $\kappa=f\circ \theta$ by uniqueness
of product integrals (cf.\ \cite[Lemma~7.4]{Mil}),
entailing that
$f(\theta([0,1]))=\kappa([0,1])$ is a compact
subset of~$(H_n)_0$.
Note that if $y\in V$, then $Q\tau(y)\cap H_0\not=\emptyset$
if and only if $\tau(y)\in H_0$, in which case
$Q\tau(y)$ is an open subset of~$H_0$.
Let
\[
D\, :=\, \{y\in V\colon Q\tau(y)\cap H_0\not=\emptyset\}\,.
\]
Then $(Q\sigma(y))_{y\in D}$
is an open cover of
$f(\theta([0,1]))$ by disjoint open subsets of~$H_0$.
Since $f(\theta([0,1]))$ is connected,
it follows that $f(\theta([0,1]))\sub Q\sigma(0)=Q$.
Hence $f(q)=f(\theta(1))\in Q$ in particular,
establishing~(\ref{ZsubQ}).\\[2.5mm]
Because $P=Q\tau(V)\isom Q\times V$
and $Q\times \{0\}$ obviously
is strongly $C^k$-initial in $Q\times V$, it follows
that $Q$ is strongly $C^k$-initial in $P$.
Since $f(Z)\sub Q$ and $f|_Z$ is $C^k$
as a map to~$P$, we infer that $f|_Z\colon Z\to Q$
is $C^k$ and hence also $f|_Z\colon Z\to H$.\\[2.5mm]
It only remains to verify the claim.
If $n=1$, there is one (and only one) choice
of $\gamma_1, \eta_1\colon \{0\}\to \{1\}$.
Now assume that $n\in \N$ and
assume that $\gamma_1,\ldots,\gamma_n$
and $\eta_1,\ldots, \eta_n$ have already
been constructed with the desired properties.
Let $A \sub R_n$ be an open $0$-neighbourhood
and $\mu \colon A \to H_{n+1}$
be a smooth map such that $\mu(0)=1$ and
$T_0(\mu) =\id_{R_n}$,
identifying $T_0(R_n)$ with $R_n$.
Consider the map
\[
\alpha\colon A\times U^n_{t_n}\times V^n_{t_n}\to G_{n+1}\,,\quad
(a,x,y)\mto \mu(a) g_n(x,y)\,.
\]
Then
\begin{equation}\label{getimag}
\im(T_{(0,x,y)}(\alpha))
\;=\; \gamma_n(x) \cdot (R_n\oplus L(G_n))\cdot \eta_n(y)
\end{equation}
for all $x\in U^n_{t_n}$ and $y\in V^n_{t_n}$
(where $\cdot$ denotes multiplication
in the tangent group $TG_{n+1}$).
In fact, abbreviating $c:=\gamma_n(x)$ and $e:=\eta_n(y)$,
the image is
\begin{eqnarray*}
R_n\cdot ce + \im \,(T_{(x,y)}g_n)
& = & R_n \cdot ce + T_{ce}G_n \; =\;
R_n\cdot ce + c \cdot L(G_n)\cdot e\\
&= & R_n \cdot ce+ c \cdot L(H_n) \cdot e + c \cdot E_n \cdot e\\
& = & R_n \cdot ce + L(H_n) \cdot ce + c \cdot E_n \cdot e\\
&= & (R_n+L(H_n))\cdot ce + c \cdot E_n \cdot e\\
& = & L(H_{n+1}) \cdot ce + c \cdot E_n \cdot e\\
&= & c\cdot L(H_{n+1}) \cdot e + c \cdot E_n \cdot e\\
&=&
c \cdot (L(H_{n+1})+E_n) \cdot e\,,
\end{eqnarray*}
where $L(H_{n+1})+E_n=R_n\oplus L(H_n)\oplus E_n=R_n\oplus L(G_n)$.\\[2.5mm]
We now equip $L(G_{n+1})$
with an inner product,
allowing us to consider\linebreak
orthogonal complements
of vector subspaces of $L(G_{n+1})$.
Setting
\[
K_y\; :=\; (\eta_n(y)^{-1}\cdot (R_n\oplus L(G_n))\cdot \eta_n(y))^\perp
\]
for $y\in V^n_{t_n}$,
we obtain a vector subbundle $K:=\bigcup_{y\in V^n_{t_n}} \{y\}\times K_y$
of the trivial bundle $V^n_{t_n}\times L(G_{n+1})\to V^n_{t^n}$.
Because $V^n_{t_n}$ is smoothly
contractible, this vector bundle
is trivial as a smooth vector bundle
(by the $C^\infty$-version of \cite[Chapter~2, Theorem~2.4]{Hir}),
and hence there exists a fibre-preserving $C^\infty$-diffeomorphism
\[
\phi\colon
V^n_{t_n}\times S_n\to K\,,
\]
where $V^n_{t_n}\times S_n\to V^n_{t_n}$
is the trivial bundle.
Let $\pi_2\colon V^n_{t_n}\times L(G_{n+1})\to L(G_{n+1})$
be the projection onto the second factor
and $\psi\colon C \to G_{n+1}$
be a $C^\infty$-map,
defined on an open $0$-neighbourhood
$C\sub L(G_{n+1})$, such that $\psi(0)=1$
and $T_0\psi=\id_{L(G_{n+1})}$.
Then $B:=(\pi_2\circ \, \phi)^{-1}(C)$ is an open
neighbourhood of $V^n_{t_n}\times \{0\}$ in
$V^n_{t_n}\times S_n$, and
\[
\nu\colon B\to G_{n+1}\,,\quad
\nu(z):=(\psi\circ \pi_2\circ \, \phi)(z)
\]
is a smooth map
such that $\nu(y,0)=1$
for all $y\in V^n_{t_n}$ and
\[
\im(T_{(y,0)}\nu)\;=\; \im T_0\nu(y,\sbull)=K_y\,.
\]
Then the smooth map
\[
\zeta \colon U^n_{t_n} \times A \times B\to G_{n+1}\,,
\quad (x,a, y,b)\mto
\mu(a)\gamma_n(x)\eta_n(y)\nu(y,b)
\]
(for $x\in U^n_{t_n}$, $a\in A$, $y\in V^n_{t_n}$
and $b\in S_n$ such that $(y,b)\in B$) 
satisfies
\begin{eqnarray*}
\im\, (T_{(x,0,y,0)}\zeta) &= &
\gamma_n(x)\cdot (R_n\oplus L(G_n))\cdot \eta_n(y)+
\gamma_n(x)\eta_n(y)\cdot K_y\\
&= &
\gamma_n(x)\eta_n(y)\cdot \big( \, \eta_n(y)^{-1} \cdot (R_n\oplus L(G_n))\cdot
\eta_n(y)\\
& & \;\;\;\; + \; (\eta_n(y)^{-1}\cdot (R_n\oplus L(G_n))\cdot \eta_n(y))^\perp
\, \big)\\
&=& \gamma_n(x)\eta_n(y)\cdot L(G_{n+1})
\;=\;
T_{\zeta(x,0,y,0)}G_{n+1}\, ,
\end{eqnarray*}
using (\ref{getimag}) for the first equality.
Hence $T_{(x,0,y,0)}\zeta$ is a linear isomorphism
(by reasons of dimension).
The restriction of $\zeta$ to $\wb{U^n_{t_{n+1}}}\times
\{0\}\times \wb{V^n_{t_{n+1}}}\times \{0\}$ corresponds to
the restriction of
$g_n$ to $\wb{U^n_{t_{n+1}}}\times \wb{V^n_{t_{n+1}}}$
and hence is injective.
Thus Lemma~\ref{tooliniti}
provides an open neighbourhood
$\Omega$ of $\wb{U^n_{t_{n+1}}}\times \{0\}\times \wb{V^n_{t_{n+1}}}\times
\{0\}$ in
$U^n_{t_n} \times A\times B$
such that $\zeta(\Omega)$ is open in $G_{n+1}$
and $\zeta|_\Omega$ is a $C^\infty$-diffeomorphism
onto $\zeta(\Omega)$.
There exists $r>0$
such that $U^n_{t_{n+1}}\times X^n_r
\times V^n_{t_{n+1}}\times Y^n_r\sub \Omega$.
Then
\[
\gamma_{n+1}\colon
U^{n+1}_{t_{n+1}}=
U^n_{t_{n+1}}\times X^n_{t_{n+1}}\to H_{n+1}\,,
\quad \gamma_{n+1}(x,a):=\mu({\textstyle\frac{r}{t_{n+1}}a})\gamma_n(x)
\]
and
\[
\eta_{n+1}\colon
V^{n+1}_{t_{n+1}}=
V^n_{t_{n+1}}\times Y^n_{t_{n+1}}\to G_{n+1}\,,
\quad \eta_{n+1}(y,b):=\eta_n(y)\nu(y,{\textstyle\frac{r}{t_{n+1}}b})
\]
are smooth maps, and
\[
g_{n+1}\colon U^{n+1}_{t_{n+1}}\times V^{n+1}_{t_{n+1}}\to G_{n+1}\,,\quad
g_{n+1}(x,y):=\gamma_{n+1}(x)\eta_{n+1}(y)
\]
is a $C^\infty$-diffeomorphism onto
an open subset of $G_{n+1}$
(because so is $\zeta|_\Omega$).
Consequently, $\gamma_{n+1}$
is a $C^\infty$-diffeomorphism
onto an open subset of $H_{n+1}$.
By construction,
$\gamma_{n+1}$ and $\eta_{n+1}$
also have all other required properties.\,\Punkt
\section{Direct limit groups are topological groups
with Lie algebra}\label{secgpLiea}
Given a (Hausdorff) topological group~$G$,
let $C(\R, G)$ be the set of
continuous $G$-valued functions on~$\R$,
and $\Hom_c(\R,G)$ be the subset of all
continuous homomorphisms from $(\R,+)$ to~$G$.
We equip $C(\R,G)$
and $\Hom_c(\R,G)$
with the compact-open topology.
Given $r\in \R$ and $\gamma \in \Hom_c(\R,G)$,
we define\linebreak
$r\gamma\in \Hom_c(\R,G)$ via $(r\gamma)(t):=\gamma(tr)$
for $t\in \R$. Then
\begin{equation}\label{defscalar}
\R\times \Hom_c(\R,G)\to \Hom_c(\R,G)\,,\quad (r,\gamma)\mto r\gamma
\end{equation}
is a continuous map~\cite[p.\,111--112]{HaM}.
We recall from \cite[Definition~2.11]{HaM}:
\begin{defn}\label{deftopliea}
$G$ is called
a \emph{topological group with Lie algebra}
if
\begin{itemize}
\item[(a)]
$(\gamma+\eta)(t):=\lim_{n\to\infty} (\gamma(t/n)\eta(t/n))^n$
exists for all $\gamma,\eta\in \Hom_c(\R,G)$,
and maps $\gamma+\eta\colon \R\to G$
so obtained are continuous homomorphisms;
\item[(b)]
$\Hom_c(\R,G)$ is a topological vector space
with the addition defined in~(a)
and the scalar multiplication (\ref{defscalar}); and
\item[(c)]
There exists a continuous bilinear map
\[
[.,.]\colon \Hom_c(\R,G)\times \Hom_c(\R,G)\to \Hom_c(\R,G)
\]
making it a Lie algebra, and such that
\[
[\gamma,\eta](t^2)\;=\;\lim_{n\to\infty} (\gamma(t/n)\eta(t/n)\gamma(-t/n)\eta(-t/n))^{n^2}\,,
\]
for all $\gamma,\eta\in \Hom_c(\R,G)$
and $t\in \R$.
\end{itemize}
\end{defn}
It is known that a Lie group~$G$ modelled
on a locally convex space is a topological group
with Lie algebra provided that~$G$ is \emph{locally
exponential} in the sense
that $G$ has a smooth exponential map $\exp_G\colon L(G)\to G$
which is a local diffeomorphism at~$0$
(see \cite[Remark~IV.1.22]{Jap}).
In particular, every Banach-Lie group
and every finite-dimensional Lie group
is a topological group with Lie algebra.
Direct limits of finite-dimensional
Lie groups always have a smooth exponential
map, but they need not be locally exponential
(see \cite[Example 5.5]{KYO}).
Nonetheless, direct limit groups are topological
groups with Lie algebra, as we verify now.
\begin{numba}\label{goodgroup}
Recall that a Lie group~$G$ is said to \emph{have an exponential map}
if, for each $x\in L(G)$,
there exists a (necessarily unique) smooth homomorphism
$\gamma_x\colon \R\to G$ with derivative
$\gamma_x'(0)=x$.
In this case, define
\mbox{$\exp_G\colon L(G)\to G$,}
$\exp_G(x):=\gamma_x(1)$.
Then $\gamma_x(t)=\exp_G(tx)$.
As a first step towards our goal,
consider a Lie group $G$ with the following
properties:
\begin{itemize}
\item[(a)]
$G$ has an exponential map;
\item[(b)]
Every continuous homomorphism $\R\to G$
is smooth;
\item[(c)]
$(\gamma_x(t/n)\gamma_y(t/n))^n\to \gamma_{x+y}(t)$ as $n\to\infty$ and
\[
(\gamma_x(t/n)\gamma_y(t/n)\gamma_x(-t/n)\gamma_y(-t/n))^{n^2}\to
\gamma_{[x,y]}(t^2)\, ,
\]
for all $x,y\in L(G)$ and $t\in \R$.
\end{itemize}
Then the map
\[
\Gamma_G \colon L(G)\to \Hom_c(\R,G)\,,\quad x\mto \gamma_x
\]
is a bijection. Furthermore, if we give $\Hom_c(\R,G)$ the unique
Lie algebra structure which makes~$\Gamma_G$
an isomorphism of Lie algebras, then conditions~(a)
and (c) from
Definition~\ref{deftopliea}
are satisfied, and the scalar multiplication
of the Lie algebra $\Hom_c(\R,G)$
is given by~(\ref{defscalar}).
Consequently, \emph{$G$ will be a topological
group with Lie algebra if $\Gamma_G$ is
a homeomorphism.}
\end{numba}
\begin{numba}\label{settDL}
Now let $G_1\sub G_2\sub\cdots$ be
finite-dimensional Lie groups,
such that the inclusion maps $G_n\to G_{n+1}$ are
smooth homomorphisms for all $n\in \N$. Let
$G=\dl\,G_n$\vspace{-.7mm} be the direct limit
Lie group.
Then (a) and (c) hold (see \cite[Proposition~4.6 (a) and (b)]{FUN}),
and also (b) (cf.\ \cite[Proposition~4.6\,(c)]{FUN}).
Moreover, the map $\Gamma_{G_n}\colon L(G_n)\to
\Hom_c(\R,G_n)$ is a homeomorphism
(cf.\linebreak
\cite[Remark~IV.1.22]{Jap}).
Since the inclusion map
$\Hom_c(\R,G_n)\to \Hom_c(\R,G)$
is continuous and $L(G)=\dl\,L(G_n)$\vspace{-.7mm}
as a topological space,
$\Gamma_G$ is continuous.
\end{numba}
\begin{thm}\label{main2}
For each direct limit Lie group
$G=\dl\,G_n$\vspace{-.7mm}
as in {\bf \ref{settDL}},
the mapping $\Gamma_G\colon L(G)\to \Hom_c(\R,G)$
is a homeomorphism
and thus $G$ is a topological group with Lie algebra.
\end{thm}
In view of the reduction steps
performed in {\bf\ref{goodgroup}}
and {\bf\ref{settDL}},
it only remains to show that $\Gamma_G$
is an open map.
This will follow from the next lemma.
\begin{numba}\label{furtred}
It is useful to recall some standard facts
first. Let $G$ be a topological group.
By definition, the sets
\[
\lfloor K, U\rfloor_G\; :=\;
\{\gamma\in C(\R,G)\colon \gamma(K)\sub U\}
\]
form a subbasis of the compact-open
topology on $C(\R,G)$, for
$K$ ranging through the compact subsets
of~$\R$ and $U$ through the open
subsets of~$G$.
It is known that
$C(\R,G)$
is a topological group under
pointwise multiplication: $(\gamma\cdot \eta)(t):=
\gamma(t)\eta(t)$.
As a consequence,
\[
\gamma \cdot \lfloor K,U\rfloor_G
\]
is an open neighbourhood of $\gamma$
for each $\gamma\in C(\R,G)$, compact set $K\sub \R$,
and identity neighbourhood $U\sub G$.
In fact, a basis of neighbourhoods of~$\gamma$
is obtained in this way.
Consequently, if $\gamma\in \Hom_c(\R,G)$,
then
\[
(\gamma \cdot \lfloor K,U\rfloor_G)\cap \Hom_c(\R,G)
\]
is an open neighbourhood of $\gamma$ in $\Hom_c(\R,G)$,
and the latter form a basis of neighbourhoods of~$\gamma$.
\end{numba}
In the next lemma,
$G_1$ and $G_2$ are finite-dimensional
Lie groups such that $G_1\sub G_2$
and the inclusion map $G_1\to G_2$ is a smooth homomorphism.
We identify $\cg_1:=L(G_1)$ with a Lie subalgebra
of $\cg_2:=L(G_2)$ and assume that
%
\begin{equation}\label{needbeysigm}
\{x\in \cg_2\colon \exp_{G_2}(\R x)\sub G_1\}\;=\; \cg_1\,.
\end{equation}
%
%
\begin{la}\label{gpLieaII}
Let $x\in \cg_1$
and $U_j\sub \cg_j$
be a relatively compact, open neighbourhood
of $x$ for $j\in \{1,2\}$,
such that
$U_1\sub U_2$.
Let $\ve>0$ and $J:=\;[-\ve,\ve]\sub \R$.
Let $V_1\sub G_1$ be a closed identity
neighbourhood such that
%
\begin{equation}\label{get3b}
(\gamma_x \cdot \lfloor J , V_1\rfloor_{G_1})\cap\Hom_c(\R,G_1)
\;\sub\;  \Gamma_{G_1}(U_1)
\end{equation}
and
$R_1 := \{y \in \cg_1\colon \exp_{G_1}(J y)\sub V_1\}$
is compact.
Then there exists
a closed identity neighbourhood
$V_2\sub G_2$ such that
\begin{equation}\label{get4}
(\gamma_x\cdot \lfloor J , V_2\rfloor_{G_2})\cap
\Hom_c(\R, G_2)
\;\sub\; \Gamma_{G_2}(U_2)\,,
\end{equation}
and $V_1$ is contained in the interior
$V_2^0$ of $V_2$ relative~$G_2$.
Furthermore, one can achieve that
\[
R_2\, := \, \{y \in \cg_2\colon \exp_{G_2}(J y)\sub V_2\}
\]
is compact.
If $V_1$ is compact, then also $V_2$
can be chosen compact.
\end{la}
\begin{proof}
Let $Q_1\supseteq Q_2 \supseteq\cdots$
be a basis of compact identity neighbourhoods in~$G_2$,
and
\[
W_m\; :=\; V_1 Q_m\quad \mbox{for $m\in \N$.}
\]
Then (\ref{get4})
holds with $V_2:=W_m$
for all sufficiently large $m\in \N$.
In fact, otherwise
we could find positive integers $m_1<m_2<\cdots$
and elements $\eta_n\in \lfloor J , W_{m_n}\rfloor$
for $n\in \N$ such that
\[
\zeta_n\, :=\, \gamma_x\cdot \eta_n\, \in \, \Hom_c(\R,G_2)
\]
and $\zeta_n\not\in \Gamma_{G_2}(U_2)$.
Then $\zeta_n=\Gamma_{G_2}(z_n)$ with
$z_n:=\Gamma_{G_2}^{-1}(\zeta_n)\in\cg_2\setminus U_2$.\\[2.5mm]
Case 1: \emph{If $(z_n)_{n\in \N}$ has a convergent
subsequence}, then we may assume (after passage to the latter)
that $z_n\to z$ as $n\to \infty$
for some $z\in \cg_2$.
Define $\eta :=\gamma_x^{-1}\cdot \gamma_z \in C(\R,G_2)$.
Then
%
\begin{equation}\label{doubleuse}
\eta(J ) \; \sub \; V_1\,.
\end{equation}
In fact, for each $t\in J$
and $n\in \N$, we have
\begin{eqnarray*}
\eta(t) & = &
\gamma_x(t)^{-1}\exp_{G_2}(tz)
\; =\; \lim_{k\to\infty}\gamma_x(t)^{-1}\exp_{G_2}(t z_k)\\
& = &\lim_{k\to\infty}\gamma_x(t)^{-1}\zeta_k(t)
\; =\; \lim_{k\to\infty}\eta_k(t)
\; \in \; W_{m_n}\,,
\end{eqnarray*}
because $\eta_k(t) \in W_{m_k}\sub W_{m_n}$ and
$W_{m_n}$ is closed.
Hence
\[
\eta(t)
\;\in\; \bigcap_{n\in \N}W_{m_n}\;=\; V_1\,,
\]
as asserted (where the last equality
holds by \cite[Lemma~3.17]{Str}).\\[2.5mm]
As a consequence of (\ref{doubleuse}),
we have $\gamma_z(t)=\gamma_x(t)\eta(t)\in G_1$
for each $t\in J$ and hence
$\exp_{G_2}(\R z)=\gamma_z(\R)\sub G_1$,
whence $z\in \cg_1$, by~(\ref{needbeysigm}).
Hence $\gamma_z \in \Hom_c(\R, G_1)$
and thus $\eta\in C(\R,G_1)$.
By~(\ref{doubleuse}),
we have $\eta\in \lfloor J, V_1\rfloor_{G_1}$.
Hence
\[
\gamma_z=\gamma_x\eta\in (\gamma_x\lfloor J , V_1\rfloor_{G_1})\cap
\Hom_c(\R, G_1)\,,
\]
whence
$\gamma_z\in \Gamma_{G_1}(U_1)$
(by (\ref{get3b}))
and thus $z\in U_1$.
But $z\in \cg_2\setminus U_2$
and hence $z\not\in U_1$, contradiction.\\[2.5mm]
Case 2:
\emph{If $(z_n)_{n\in \N}$ has no convergent
subsequence}, after passing to a subsequence
we may assume that $\|z_n\|\to\infty$ as $n\to\infty$,
where $\|.\|$ is a given\linebreak
norm on~$\cg_2$.
Since $R_1$ is compact,
there exists $r>0$ such that
%
\begin{equation}\label{hhelp}
R_1\; \sub \; \{y\in \cg_1\colon \|y\|<r\}\; =:\; B_1\, .
\end{equation}
We let $B_2:= \{y\in \cg_2\colon \|y\|<2r\}$.
As $2r z_n/\|z_n\|$
is contained in
the~compact set $\partial B_2$,
after passage to a subsequence
we may assume that
$2r z_n/\|z_n\|\to b$ as $n\to\infty$
for some $b\in \partial B_2$.
Abbreviate $r_n:=2r/\|z_n\|$.
Then $r_n\to 0$, and
we may assume that $r_n\leq 1$ for each~$n$.
For each $t\in J$,
we then have
\[
\eta_n(t r_n)\;=\;
\exp_{G_2}(-tr_n x)\exp_{G_2}(tr_nz_n)\;\to\;
\exp_{G_2}(tb)\;\;  \mbox{as $\; n\to\infty$.}
\]
Since $\eta_n(tr_n)\in W_{m_n}$,
it follows that
$\exp_{G_2}(tb)\in \bigcap_{n\in \N}W_{m_n}=V_1$.
Hence
$\exp_{G_2}(\R b)\sub G_1$
and thus $b\in \cg_1$.
By the preceding, we have $b\in R_1$
and therefore $b\in B_1$, by~(\ref{hhelp}).
But $b\in \partial B_2$,
whence $b\not\in B_2$ and
thus $b\not\in B_1$,
contradiction.\\[2.5mm]
The final assertion is clear,
since compactness of~$V_1$
entails that also each~$W_m$ is compact.
It only remains to show that~$R_2$
can be chosen compact.\linebreak
However, applying the results
already shown with $1$, $B_1$ and $B_2$
instead of~$x$, $U_1$ and $U_2$, respectively,
we see that
\[
\{y\in \cg_2\colon \exp_{G_2}(J y)\sub W_m\}
\sub \; B_2
\]
for all sufficiently large~$m$.
For such~$m$,
the set on the left
is compact.
\end{proof}
{\bf Proof of Theorem~\ref{main2}.}
Since $G_0=\dl\,(G_n)_0$\vspace{-.8mm} (see {\bf\ref{fct1}}),
$L(G)=L(G_0)$,
$\Hom_c(\R,G)=\Hom_c(\R, G_0)$
and $\Gamma_G=\Gamma_{G_0}$,
after replacing
$G$ with $G_0$ and each $G_n$ with
$G_{n,0}$ we may assume that $G$ and each $G_n$
is connected. The discussion leading to (\ref{givesgoodH})
now shows that
\[
L(G_n)\; =\; \{x\in L(G_{n+1})\colon \exp_{G_{n+1}}(\R x)\sub G_n\}
\]
(as needed in Lemma~\ref{gpLieaII}).\\[2.5mm]
It remains to show that $\Gamma_G(U)$
is open in $\Hom_c(\R,G)$
for each open subset
$U\sub L(G)=\bigcup_{n\in \N}L(G_n)$.
We verify that $\Gamma_G(U)$
is a neighbourhood of $\Gamma_G(x)$
for each $x\in U$.
After passage to a cofinal subsequence,
we may assume that $x\in L(G_1)$.
Then $U_n:=L(G_n)\cap U$ is an
open neighbourhood of~$x$ in~$L(G_n)$
for each $n\in \N$, and $U_1\sub U_2\sub \cdots$.
Since $\Gamma_{G_1}(U_1)$ is an open neighbourhood of $\gamma_x$
in $\Hom_c(\R,G_1)$, there exists $\ve>0$ and an open identity
neighbourhood $P\sub G_1$
such that
\[
(\gamma_x\lfloor J, P\rfloor_{G_1})\cap \Hom_c(\R, G_1)\; \sub\;
\Gamma_{G_1}(U_1)\,,
\]
where $J:=[{-\ve},\ve]$ (cf.\ {\bf\ref{furtred}}).
There exists a compact identity
neighbourhood $V_1\sub P$
such that
\[
R_1\; :=\; \{y\in L(G_1)\colon \exp_{G_1}(J y)\sub V_1\}
\]
is compact.
Using now Lemma~\ref{gpLieaII}
and induction, we find compact identity
neighbourhoods $V_n\sub G_n$ for all integers
$n\geq 2$ such that
\begin{itemize}
\item[(i)]
$V_{n-1}\sub V_n^0$ (the interior relative $G_n$);
\item[(ii)]
$(\gamma_x\cdot \lfloor J, V_n\rfloor_{G_n})\cap \Hom_c(\R, G_n)
\sub \Gamma_{G_n}(U_n)$, and
\item[(iii)]
$R_n:=\{y\in L(G_n)\colon \exp_{G_n}(J y)\sub V_n\}$
is compact.
\end{itemize}
Since $Q_n:=V_n^0$ is open in~$G_n$
and $Q_1\sub Q_2\sub\cdots$, the set
$Q:=\bigcup_{n\in \N}Q_n$ is open
in $G$ (see {\bf \ref{fct2}}).
Then
\[
(\gamma_x\cdot \lfloor J, Q\rfloor_G)\cap \Hom_c(\R, G)
\]
is an open neighbourhood of $\gamma_x$ in $\Hom_c(\R, G)$.
The proof will be complete
if we can show that
\begin{equation}\label{willend}
(\gamma_x\cdot \lfloor J, Q\rfloor_G)\cap \Hom_c(\R, G)\;\sub\;
\Gamma_G(U)\,.
\end{equation}
To verify (\ref{willend}), let
$\eta\in \lfloor J,Q\rfloor_G$ such that
$\zeta:=\gamma_x\cdot\eta\in \Hom_c(\R,G)$.
Then
$\zeta=\gamma_z$ with $z:=\Gamma_G^{-1}(\zeta)$.
There is $n\in \N$ such that
$z\in L(G_n)$ and thus $\zeta\in \Hom_c(\R,G_n)$.
Then $\eta=\gamma_x^{-1}\cdot \zeta \in C(\R, G_n)$.
Since $\eta(J)\sub G_n$ is compact
and $(Q_m\cap G_n)_{m\geq n}$ is an open cover of
$\eta(J)$, there exists $m\geq n$ such that
$\eta(J)\sub Q_m$.
Thus $\eta\in \lfloor J,Q_m\rfloor_{G_m}\sub
\lfloor J, V_m\rfloor_{G_m}$
and hence
\[
\zeta\, = \, \gamma_x\cdot \eta\, \in\,
(\gamma_x \cdot \lfloor J,V_m\rfloor_{G_m})\cap \Hom_c(\R,G_m)\, \sub \, \Gamma_{G_m}(U_m)\,,
\]
using (ii).
Then $\zeta\in \Gamma_{G_m}(U_m)\sub \Gamma_G(U)$
and thus (\ref{willend}) is established.\,\vspace{2mm}\Punkt

\noindent
If a Lie group is a topological group
with Lie algebra, then this has
useful consequences.
We recall that it is an unsolved open
problem (first formulated by John Milnor)
whether every continuous
homomorphism between infinite-dimensional
Lie groups is smooth (as in the finite-dimensional
case). The following result provides
some positive information.
\begin{prop}
Let $H$ be a locally exponential Lie group
and $G$ be a Lie group such that
{\rm (a)--(c)} from {\bf\ref{goodgroup}}
are satisfied, $\Gamma_G\colon L(G)\to \Hom_c(\R,G)$
is a homeomorphism, and $\exp_G\colon L(G)\to G$
is smooth.
Then every continuous homomorphism
$\alpha\colon H\to G$ is smooth.
\end{prop}
\begin{proof}
The map
$\beta := \Hom_c(\R,\alpha)\colon
\Hom_c(\R,H)\to\Hom_c(\R,G)$,
$\gamma\mto\alpha\circ \gamma$
is continuous.
Since
\begin{eqnarray*}
\beta(\gamma+\eta)(t)
&= & \alpha\left(
\lim_{n\to\infty}(\gamma(t/n)\eta(t/n))^n\right)\\
&= &
\lim_{n\to\infty} ((\alpha\circ \gamma)(t/n)(\alpha\circ \eta)(t/n))^n
\; =\; (\beta(\gamma) +\beta(\eta))(t)
\end{eqnarray*}
for all $\gamma,\eta\in \Hom_c(\R,H)$,
we see that $\beta$ is a homomorphism
of groups and hence
a continuous linear map.
Then $\theta:=\Gamma_G^{-1}\circ \beta\circ \Gamma_H\colon
L(H)\to L(G)$ is a continuous linear (and hence smooth)
map such that $\exp_G\circ \, \theta=\alpha\circ \exp_H$.
Since $\exp_G\circ \, \theta$ is smooth
and $\exp_H$ is a local diffeomorphism, it
follows that~$\alpha$ is smooth on some
open identity neighbourhood and hence smooth.
\end{proof}
\section{Borel-type theorem for diffeomorphisms}\label{secborel}
This section is devoted to the proof of Borel's
Theorem for Diffeomorphisms,
as stated in the introduction.\\[2.5mm]
As usual, $\GL_n(\R)$ denotes the group
of all invertible $n \!\times\! n$-matrices
and $\GL(\R^n)$ the group of all automorphisms
of the real vector space~$\R^n$.
The algebra of all linear endomorphism
of~$\R^n$ will be denoted by $\cL(\R^n)$.
If~$M$ is a finite-dimensional
smooth manifold and $K\sub M$,
we let $\Diff_K(M)$ be the Lie group
of all $C^\infty$-diffeomorphisms $\gamma\colon M\to M$
such that $\gamma(x)=x$ for all $x\in M\setminus K$
(see, e.g., \cite{DFM}; cf.\ \cite{Mic}).
We write $C^\infty_K(M)$ for the space of smooth
maps $\gamma\colon M\to\R$ such that
$\gamma|_{M\setminus K}=0$.
Finally, we let $C^\infty(\R^n,\R^n)$
be the space of $\R^n$-valued smooth
maps on~$\R^n$, equipped with the
usual locally convex topology
(the smooth compact-open topology).
\begin{la}
Let $n\in \N$, $A\in \GL(\R^n)$
such that $\det(A)>0$, and $K\sub \R^n$
be a compact $0$-neighbourhood.
Then there exists a smooth map
\[
\zeta\colon \R\to \Diff_K(\R^n)
\]
with $\zeta(0)=\id_{\R^n}$, such that
\[
\theta:=\zeta(1)
\]
satisfies $\theta|_U=A|_U$ for some
$0$-neighbourhood $U\sub K$.
\end{la}
\begin{proof}
Let $R$ be a compact
$0$-neighbourhood
contained in the interior $K^0$ of~$K$.
Then there is $\chi\in C^\infty_K(\R^n)$
such that $\chi|_R=1$.
It is well known that
\[
\{B\in \GL(\R^n)\colon \det(B)>0\}
\]
is the connected component of
$\GL(\R^n)$. We therefore
find a smooth map $f\colon \R\to \GL(\R^n)$
with
$f(0)=\id_{\R^n}$
and $f(1)=A$.
Let $f'(t)\in \cL(\R^n)$
be the derivative of~$f$ at $t\in \R$.
We consider the time-dependent smooth vector
field
\[
v\colon \R\times \R^n\to \R^n\,,\quad
v(t,x):=(f'(t)\circ f(t)^{-1})(x)
\]
on $\R^n$. For each $x\in \R^n$,
\[
\gamma_x\colon \R\to \R^n\,,\quad \gamma_x(t):=f(t)(x)
\]
is a smooth curve such that $\gamma_x(0)=x$
and
\[
\gamma_x'(t)=f'(t)(x)=(f'(t)\circ f(t)^{-1})(f(t)(x))
=v(t,\gamma_x(t))\,.
\]
Thus $\gamma_x$ solves the initial value problem
\[
y'(t)=v(t,y(t))\,,\quad y(0)=x\,.
\]
Hence $\Phi_v\colon \R\times \R^n\to \R^n$, $\Phi_{v,t}(x):=\gamma_x(t)$
is the flow of~$v$ (for fixed initial time $t_0=0$).
Now consider
\[
w\colon \R\times \R^n\to \R^n\,,\quad
w(t,x):=\chi(x)v(t,x)\,.
\]
Then $w$ admits a smooth flow $\Phi_w\colon \R\times \R^n\to \R^n$,
and $\zeta(t):=\Phi_{w,t}\in \Diff_K(\R^n)$ for each $t\in \R$.
Since $\phi_{w,t}(0)=0$ for each $t\in \R$
and the flow is continuous, we find a
$0$-neighbourhood $U\sub K$ such that
\[
\Phi_{w,t}(x)\in R^0\quad\mbox{for all $x\in U$ and $t\in [0,1]$.}
\]
Since $v$ and $w$ coincide on $\R\times R^0$,
it follows that $\Phi_{w,t}(x)=\Phi_{v,t}(x)$
for all $x\in U$ and $t\in [0,1]$,
entailing that
\[
\Phi_{w,1}(x)=\Phi_{v,1}(x)=\gamma_x(1)=A(x)\quad\mbox{for all $x\in U$.}
\]
By construction of the Lie group structure
on $\Diff_K(\R^n)$,
the map
\[
\kappa\colon \Diff_K(\R^n)\to C^\infty_K(\R^n,\R^n)\, ,\quad
\gamma\mto \gamma-\id_{\R^n}
\]
is a diffeomorphism onto an open subset
of $C^\infty_K(\R^n,\R^n)$ (cf.\ \cite{DRN}).
Because
\[
\R\times \R^n\to \R^n\, , \quad (t,x)\mto (\zeta(t)-\id_{\R^n})(x)
=\Phi_{w,t}(x)-x
\]
is a smooth mapping, the exponential law for smooth mappings
(see, e.g., \cite[Lemma~12.1\,(a)]{ZOO}) entails that $\kappa\circ \zeta$
(and hence also $\zeta$) is smooth.
Since $\zeta(0)=\id_{\R^n}$
and $\zeta(1)=\Phi_{w,1}$,
we see that $\zeta$ has all required properties.
\end{proof}
Given $n\in \N$ and $j\in \{1,\ldots, n\}$,
let $e_j:=(0,\ldots, 0, 1,0,\ldots, 0)^T\in \R^n$
with $1$ in the $j$-th slot.
\begin{la}\label{seclemdiff}
Let $n\in \N$
and $a_\alpha\in \R^n$
for multi-indices $\alpha\in \N_0^n$
be given such that
$A:=(a_{e_1},\ldots, a_{e_n})\in \GL_n(\R)$
and $\det(A)>0$.
Then there exists $\phi\in \Diff_c(\R^n)$
such that $\frac{(\partial^\alpha\phi)(0)}{\alpha!}=a_\alpha$
for each $\alpha\in \N_0^n$.
\end{la}
\begin{proof}
After a translation, we may assume that $a_0=0$.
Then the vectors $a_\alpha$
determine a formal diffeomorphism $f=\sum_\alpha a_\alpha x^\alpha$ of $\R^n$.
Considering~$A$ as a formal diffeomorphism,
the formal composition $g:=A^{-1}\circ f$ is of the form
$g=\sum_\alpha b_\alpha x^\alpha$
with vectors $b_\alpha\in \R^n$
such that $(b_{e_1},\ldots, b_{e_n})={\bf 1}$
is the identity matrix.
Let $\theta$ be as in the previous lemma.
If we can prove the lemma with the $b_\alpha$
in place of the $a_\alpha$,
leading to~$\phi$,
then $\kappa:=\theta \circ \phi$
is a diffeomorphism
such that
$\kappa|_W= A\circ \phi|_W$
for some $0$-neighbourhood $W\sub \R^n$,
and hence
$\frac{(\partial^\alpha\kappa)(0)}{\alpha!}=a_\alpha$
for each $\alpha\in\N_0^n$.
We may therefore assume now
that $A=\one$.\\[2.5mm]
Let $\|.\|$ be the maximum norm on $\R^n$.
Then $\{x\in \R^n\colon \|x\|\leq \ve\}\sub K$
for some $\ve>0$.
Let $h\colon \R\to \R$
be a compactly supported smooth
function such that $h(x)=x$
for $x$ in some $0$-neighbourhood,
and $h(x)=0$ if $|x| \geq \ve$.
Define $h_\alpha(x_1,\ldots, x_n):=h(x_1)^{\alpha_1}\cdots
h(x_n)^{\alpha_n}$
and
\[
M_{k,m}\; :=\; \max_{|\beta|=m} \, \sum_{|\alpha|=k} \, \|a_\alpha\|
\cdot \|\partial^\beta h_\alpha\|_\infty\,,
\]
where $\|.\|_\infty$ is the supremum norm
on the space $C_0(\R)$ real-valued functions
that vanish at infinity,
and $|\alpha|=\alpha_1+\cdots +\alpha_n$
the order of the multi-index $\alpha\in\N_0^n$.
For each $k\in \N$
with $k\geq 2$,
pick $c_k>1$ such that
\[
c^{m-k}_k M_{k,m}\;<\; \frac{2^{-k}}{n}
\]
for all $m\in \N_0$
such that $m<k$.
Define $\phi_k\colon \R^n\to\R^n$ via
$\phi_k(x):=\sum_{|\alpha|=k}a_\alpha c_k^{-k}h_\alpha(c_kx)$.
Then $\phi_k$ is smooth
and
\[
\partial^\beta\phi_k(x)
\;=\; c_k^{|\beta|-k} \sum_{|\alpha|=k}a_\alpha 
(\partial^\beta h_\alpha)(c_kx)\,,
\]
whence
\[
\|\partial^\beta\phi_k\|_\infty
\; \leq \; c_k^{|\beta|-k}M_{k,|\beta|}\, ,
\]
which is $\leq \frac{2^{-k}}{n}$
for all $k\in \N$ such that $k>|\beta|$.
Hence the limit
$\psi:=\sum_{k=2}^\infty \phi_k$
exists in $C^\infty(\R^n,\R^n)$,
and clearly $\psi$ has compact support.
By the preceding, for $|\beta|=1$
we have
$\|\partial^\beta\phi_k\|_\infty
\leq \frac{2^{-k}}{n}$
for each $k\geq 2$
and thus
$\|\partial^\beta \psi\|_\infty\leq\sum_{k=2}^\infty\frac{2^{-k}}{n}
=\frac{1}{2n}$,
whence
$\|\psi'(x)\|_{\text{op}}\leq 1/2$
for each $x\in \R^n$
and thus
$\sup\{\|\psi'(x)\|_{\text{op}}\colon x\in \R^n\}
\leq 1/2<1$.
Therefore
$\phi:=\id_{\R^n}+\psi \in \Diff_c(\R^n)$
(see \cite[Lemma~5.1]{DRN}).
By construction,
$\frac{\partial^\alpha\phi(0)}{\alpha!}=a_\alpha$
for each $\alpha$.
\end{proof}
{\bf Proof of Borel's Theorem for Diffeomorphisms.}
Let
$p$ and a chart $\phi\colon
M\supseteq U\to V\sub \R^n$
be as described in the theorem.
After shrinking~$U$, we may assume that~$V$
is a euclidean ball
with center~$0$. Composing now with a diffeomorphism
$V\to \R^n$ which is the identity
on some $0$-neighbourhood and fixes~$0$,
we can (and shall) assume instead
that $V=\R^n$.
Let $C\sub U$ be a compact neighbourhood of~$p$
and $K := \phi(C)$.
Then the restriction map
\[
\Diff_C(M)\to \Diff_C(U)
\]
is an isomorphism of Lie groups
and so is the map
\[
\Diff_C(M)\to \Diff_K(\R^n)\,,\quad \gamma\mto \phi\circ \gamma\circ \phi^{-1}
\]
(this
is clear from the construction of the Lie group
structure in \cite{DFM}).
It therefore suffices to consider
the case where $M=\R^n$, $\phi=\id_{\R^n}$
and $p=0$.
In this case, the surjectivity of
\[
T^\infty_p\colon \Diff_c(\R^n)_{p,0}\to \Gf_n(\R)_0
\]
has been established
in Lemma~\ref{seclemdiff}.
The Chain Rule for Taylor Polynomials
entails that $T^\infty_p$ is a homomorphism
of groups. Since $\Gf_n(\R)_0$ is an open subset
of $\{\sum_\alpha a_\alpha x^\alpha\in
\R[\![x_1,\ldots, x_n]\!]\colon a_0=0\}
\isom \R^{\N_0^n\setminus \{0\}}$
with the product topology,
and each component
\[
(T^\infty_p)_\alpha\colon \Diff_c(\R^n)_{p,0}\to\R\, ,\quad
\gamma\mto \frac{(\partial^\alpha\gamma)(0)}{\alpha!}
\]
is smooth, we deduce that $T^\infty_p$ is smooth.
This completes the proof.\,\vspace{3mm}\Punkt

\noindent
\emph{Acknowledgements}.
The author thanks K.-H. Neeb for
communicating early versions of his survey,
and for comments on a draft
of the current article.
\noindent
{\footnotesize
{\bf Helge Gl\"{o}ckner}, Universit\"{a}t Paderborn,
Institut f\"{u}r Mathematik,
Warburger Str.\ 100,\\ 33098 Paderborn,
Germany.
\,E-Mail: {\tt glockner@math.uni-paderborn.de}}

\begin{thebibliography}{99}
%
%
\bibitem{Ber}
Bertram, W., H. Gl\"{o}ckner
and K.-H. Neeb,
\emph{Differential calculus over general
base fields and rings}, Expo.\ Math.\ {\bf 22} (2004),
213--282.
%
%
\bibitem{Bou}
Bourbaki, N., ``Lie Groups and Lie Algebras, Chapters 1--3,''
Springer-Verlag, 1989.
%
%
\bibitem{DaG} Dahmen, R. and H. Gl\"{o}ckner,
\emph{Regularity in Milnor's sense for
direct limits of infinite-dimensional Lie groups},
in preparation.
%
%
\bibitem{RES}
Gl\"{o}ckner, H.,
\emph{Infinite-dimensional Lie groups without completeness restrictions},
pp.\,43--59 in:
A. Strasburger et al.\ (eds.)
``Geometry and analysis on finite- and infinite-dimensional Lie groups,''
Banach Center Publ.\ {\bf 55} (2002), 43--59.
%
%
\bibitem{KYO} Gl\"{o}ckner, H.,
\emph{Direct limit Lie groups and manifolds},
J. Math.\ Kyoto Univ.\ {\bf 43} (2003), 1--26.
%
%
\bibitem{DRN}
Gl\"{o}ckner, H., \emph{$\Diff(\R^n)$ as
a Milnor-Lie group},
Math.\ Nachr.\ {\bf 278} (2005),
1025--1032.
%
%
\bibitem{FUN} Gl\"{o}ckner, H.,
\emph{Fundamentals of direct limit Lie theory},
Compos.\ Math.\ {\bf 141} (2005), 1551--1577.
%
%
%
%
\bibitem{JFA}
Gl\"{o}ckner, H.,
\emph{Direct limits of infinite-dimensional Lie groups
compared to direct limits in related categories},
J. Funct.\ Anal.\  {\bf 245} (2007), 19--61.
%
%
\bibitem{DFM}
Gl\"{o}ckner, H., \emph{Patched locally convex spaces,
almost local mappings and diffeomorphism groups
of non-compact manifolds}, manuscript, 2002.
%
%
\bibitem{ZOO}
Gl\"{o}ckner, H., \emph{Lie groups
over non-discrete topological fields},
preprint, arXiv:math/0408008v1.
%
%
\bibitem{GaL}
Gl\"{o}ckner, H. and L.\,R. Lovas, \emph{Frobenius
and Stefan-Sussmann theorems for
various types of distributions on infinite-dimensional
manifolds}, manuscript in preparation.
%
%
\bibitem{GaN}
Gl\"{o}ckner, H. and K.-H. Neeb,
``Infinite-Dimensional Lie Groups, Vol.\,I'' 
book in preparation.
%
%
%
%
\bibitem{Hil} 
Hiltunen, S.,
\emph{A Frobenius theorem for locally convex global analysis},
Monatsh.\ Math.\  {\bf 129} (2000), 109--117.
%
%
\bibitem{HST}
Hirai, T., H. Shimomura, N. Tatsuuma and E. Hirai,
\emph{Inductive limits of topologies, their direct product,
and problems related to algebraic structures},
J. Math.\ Kyoto Univ.\ {\bf 41} (2001),
475--505.
%
%
\bibitem{Hir}
Hirsch, M.\,W., ``Differential Topology,''
Springer-Verlag, 1976.
%
%
\bibitem{HaM}
Hofmann, K.\,H. and S.\,A. Morris,
``The Structure of Connected Pro-Lie Groups,''
EMS Tracts in Math.\ {\bf 2},
EMS Publ.\ House,
Zurich, 2007.
%
%
\bibitem{KaM}
Kriegl, A. and P.\,W. Michor,
``The Convenient Setting of
Global Ana\-lysis,'' Amer.\ Math.\ Soc.,
Providence, 1997.
%
%
\bibitem{Mic}
Michor, P.\,W.,
``Manifolds of differentiable mappings,'' 
Shiva Publishing, Nantwich, 1980.
%
%
\bibitem{Mil}
Milnor, J.,
\emph{Remarks
on infinite-dimensional Lie groups}.
pp.\,1007--1057 in:
Relativity, groups and topology, II (Les Houches, 1983),
North-Holland, Amsterdam, 1984. 
%
%
\bibitem{NRW}
Natarajan, L., E. Rodr\'{\i}guez-Carrington
and J.\,A. Wolf,
\emph{Differentiable structure for direct limit groups},
Letters in Math.\ Phys.\ {\bf 23} (1991), 99--109.
%
%
\bibitem{Jap}
Neeb, K.-H.,
\emph{Towards a Lie theory of locally convex groups},
\mbox{Jpn.\ J. Math.\ {\bf 1}} (2006), 291--468.
%
%
\bibitem{Pis}
Pisanelli, D.,
\emph{An example of an infinite Lie group},
Proc.\ Amer.\ Math.\ Soc.\
{\bf 62} (1977), 156--160.
%
%
\bibitem{RaK}
Robart, T. and N. Kamran,
\emph{Sur la th\'{e}orie locale des pseudogroupes de transformations
continus infinis} I,
Math.\ Ann.\  {\bf 308}  (1997),  593--613.
%
%
\bibitem{Str}
Stroppel, M.,
``Locally Compact Groups,''
EMS Publ.\ House,
Zurich, 2006.\vspace{.3mm}
%
%
\end{thebibliography}
\end{document}